\documentstyle[leqno]{article}

\newtheorem{th}{Theorem}[subsection]
\newtheorem{prop}[th]{Proposition}

\newcounter{defin}[subsection]
\renewcommand{\thedefin}{\thesubsection.\arabic{defin}}
\newcounter{ex}[subsection]
\renewcommand{\theex}{\thesubsection.\arabic{ex}}
\newcounter{rem}[subsection]
\renewcommand{\therem}{\thesubsection.\arabic{rem}}
\title{The Geometry of Relativistic Rheonomic\\
Lagrange Spaces}
\author{Mircea Neagu}
\date{}
\begin{document}
\maketitle
\begin{abstract}
In this paper we shall present a geometrization of time-dependent Lagrangians.
The reader is invited  to compare this geometrization with that contained in
the book of Miron and Anastasiei \cite{11}.
In order to develope the subsequent {\it Relativistic Rheonomic Lagrange Geometry,}
Section 1 describes the main geometrical aspects of the 1-jet space $J^1(R,M)$,
in the sense  of d-tensors, d-connections, d-torsions and d-curvatures. Section 2
introduces the notion of {\it Relativistic Rheonomic Lagrange Space},  which
naturally generalizes that of {\it Classical Rheonomic Lagrange Space}
\cite{11}, and constructs its canonical nonlinear connection $\Gamma$
as well as its Cartan canonical $\Gamma$-linear connection. We point out that
our geometry gives a model for both gravitational and electromagnetic field.
From this point of view, Section 4 presents the Maxwell equations of the
relativistic rheonomic Lagrangian electromagnetism. Section 5 describes the
Einstein's gravitational field equations of a relativistic rheonomic Lagrange
space.
\end{abstract}

\section{The geometry of $J^1(R,M)$}

\subsection{Some physical aspects}

\hspace{5mm} Let us consider the usual time axis represented by the set of
real numbers $R$ and a real, smooth and $n$-dimensional manifold $M$ that
we regard like a {\it "spatial"} manifold \cite{16}. We suppose that the temporal
manifold $R$ is coordinated by $t$ while the spatial manifold $M$ is coordinated
by $(x^i)_{i=\overline{1,n}}$. Note that, throughout this paper, the latin
letters $i,j,k\ldots$ run from 1 to $n$.

Let $J^1(R,M)\equiv R\times TM$ be the usual 1-jet vector bundle, coordinated
by $(t,x^i,y^i)$, and regarded over the product manifold base $R\times M$.
From physical point of view, the fibre bundle
\begin{equation}
J^1(R,M)\to R\times M,\quad (t,x^i,y^i)\to(t,x^i),
\end{equation}
is regarded like a {\it bundle of configurations}, in mechanics terms.
The gauge group of this bundle of configurations is
\begin{equation}\label{G_1}
\left\{\begin{array}{l}
\tilde t=\tilde t(t)\\
\tilde x^i=\tilde x^i(x^j)\\
\displaystyle{\tilde y^i={\partial\tilde x^i\over\partial x^j}{dt\over
d\tilde t}y^j.}
\end{array}\right.
\end{equation}
We  remark that the form of this gauge group stands out by the {\it relativistic}
character of the time $t$. For that reason, we consider that the jet vector
bundle of order one $J^1(R,M)$ is a natural house of the {\it relativistic
rheonomic Lagrangian mechanics}.

It is important to note that, in the {\it classical rheonomic Lagrangian
mechanics} \cite{11}, the bundle of configuration is the fibre bundle
\begin{equation}
\pi: R\times TM\to M,\;(t,x^i,y^i)\to (x^i),
\end{equation}
whose geometrical invariance group is
\begin{equation}\label{G_2}
\left\{\begin{array}{l}
\tilde t=t\\
\tilde x^i=\tilde x^i(x^j)\\
\displaystyle{\tilde y^i={\partial\tilde x^i\over\partial x^j}y^j.}
\end{array}\right.
\end{equation}
The structure of the gauge group \ref{G_2} emphasizes the  {\it
absolute}  character of the time $t$ from the classical rheonomic Lagrangian
mechanics. At the same time, we point out that the gauge group \ref{G_2} is
a subgroup of \ref{G_1}. In other words, the gauge group of the jet bundle of
order one from the relativistic rheonomic Lagrangian mechanics is more general
than that used in the classical rheonomic Lagrangian mechanics, which ignores
the temporal reparametrizations.

Finally, we point out that a deeply exposition of the physical aspects of
the classical rheonomic Lagrange geometry is done by Ikeda in \cite{6} and
\cite{12}. At the same time, we invite the reader to compare the classical
rheonomic Lagrangian mechanics \cite{10}  with that relativistic, whose
geometrical background is developed in this paper.

\subsection{Time-dependent sprays. Harmonic curves}

\setcounter{equation}{0}
\hspace{5mm} Let us consider that the temporal manifold $R$ is endowed with a
semi-Riemannian metric $h=(h_{11}(t))$. In order to develope the
geometrical background of the relativistic rheonomic mechanics on the
1-jet fibre bundle $E=J^1(R,M)$, we will introduce a collection of important
geometrical concepts. An important geometrical concept on $J^1(R,M)$ is that of
time-dependent spray, which naturally generalizes the notion of time-dependent
spray on $R\times M$, used in \cite{11} and \cite{22}. In order to
introduce this concept, let us consider the following notions:\medskip\\
\addtocounter{defin}{1}
{\bf Definition \thedefin} A global tensor $H$ (resp. $G$) on $E$, locally
expressed by
\begin{equation}\label{ts}
H=dt\otimes{\partial\over\partial t}-2H^{(j)}_{(1)1}dt\otimes{\partial\over
\partial y^j},
\end{equation}
respectively
\begin{equation}\label{ss}
G=y^jdt\otimes{\partial\over\partial x^j}-2G^{(j)}_{(1)1}dt\otimes{\partial
\over\partial y^j}
\end{equation}
is called a {\it temporal} (resp. {\it spatial}) {\it spray} on $E$.\medskip

Because the sprays $H$ and $G$ are global tensors, using the coordinate
transformations \ref{G_1} on the 1-jet space $E$, it is easy to deduce the
following \cite{16}

\begin{th}
To give a temporal (spatial) spray on $E$ is equivalent to give a set of
local functions $H=(H^{(j)}_{(1)1})$ (resp. $G=(G^{(j)}_{(1)1})$) which
transform by the rules
\begin{equation}\label{lts}
2\tilde H^{(k)}_{(1)1}=2H^{(j)}_{(1)1}\left({dt\over d\tilde t}\right)^2{\partial
\tilde x^k\over\partial x^j}-{dt\over d\tilde t}{\partial\tilde y^k\over
\partial t},
\end{equation}
respectively
\begin{equation}\label{lss}
2\tilde G^{(k)}_{(1)1}=2G^{(j)}_{(1)1}\left({dt\over d\tilde t}\right)^2{\partial
\tilde x^k\over\partial x^j}-{\partial x^i\over\partial\tilde x^j}{\partial
\tilde y^k\over\partial x^i}\tilde y^j.
\end{equation}
\end{th}

The previous theorem allows us to offer the following important examples of
temporal and spatial sprays. The importance of these sprays comes from their
using in the description of the local equations of harmonic maps between two
semi-Riemannian manifolds \cite{4}.\medskip\\
\addtocounter{ex}{1}
{\bf Example \theex} Let $h=(h_{11})$  (resp. $\varphi=(\varphi_{ij})$) be a
semi-Riemannian metric on $R$ (resp. $M$) and $H^1_{11}$ (resp.
$\gamma^i_{jk}$) its Christoffel symbols. In this context, taking into account
the transformation rules of the Christoffel symbols $H^1_{11}$ and $\gamma^i_
{jk}$, we deduce that the components
$2H^{(j)}_{(1)1}=-H^1_{11}y^j$ (resp. $2G^{(j)}_{(1)1}=\gamma^j_{kl}y^ky^l$)
represent a temporal (resp. spatial) spray which is called the {\it canonical
temporal} (resp. {\it spatial}) {\it spray associated to the metric $h$}
(resp. $\varphi$).\medskip\\
\addtocounter{defin}{1}
{\bf Definition \thedefin} A pair $(H,G)$, which consists of a temporal  spray
and a  spatial one, is called a {\it time-dependent spray on $J^1(R,M)$.}
\medskip

Follwing the geometrical development of the classical rheonomic Lagrange
mechanics, we introduce a natural  generalization of the notion of path of
a time-dependent spray, used in \cite{11}.\medskip\\
\addtocounter{defin}{1}
{\bf Definition \thedefin} A curve $c\in C^\infty(R,M)$ is called a {\it
harmonic curve of the time-dependent spray $(H,G)$ on $J^1(R,M)$, with respect
to the semi-Riemannian  temporal metric $h=(h_{11}(t))$ on R}, if $c$ is a
solution of the DEs system of order two
\begin{equation}\label{hc}
h^{11}\left\{{d^2x^i\over dt^2} +2G^{(i)}_{(1)1}+2H^{(i)}_{(1)1}\right\}=0,
\end{equation}
where $h^{11}h_{11}=1$ and the curve $c$ is locally expressed by $R\ni t\to
(x^i(t))_{i=\overline{1,n}}\in M$.\medskip\\
\addtocounter{rem}{1}
{\bf Remarks \therem} i) Under the coordinate transformations of $J^1(R,M)$, the
left term of the equations \ref{hc} modifies like a d-tensor, that is,
\begin{equation}\label{thc}
\hspace*{5mm}\left[h^{11}\left\{{d^2x^i\over dt^2} +2G^{(i)}_{(1)1}+2H^{(i)}_{(1)1}\right\}
\right]=
{\partial x^i\over\partial\tilde x^j}\left[
\tilde h^{11}\left\{{d^2\tilde x^j\over d\tilde t^2} +2\tilde G^{(j)}_{(1)1}+2\tilde
H^{(j)}_{(1)1}\right\}\right].
\end{equation}
Consequently, the equations \ref{hc} are global on $J^1(R,M)\equiv R\times TM$
(i. e. their geometrical invariance group is \ref{G_1}).

ii) Comparatively, the equations of a path on $R\times TM$ (see \cite{11}),
that we generalized by \ref{hc}, are invariant only under the gauge group \ref{G_2}.
\medskip\\
\addtocounter{ex}{1}
{\bf Example \theex} Let us consider the canonical sprays asociated to the metrics
$h$ and $\varphi$, which are locally expressed by
\begin{equation}\label{cs}
\left\{\begin{array}{l}\medskip
\displaystyle{H^{(i)}_{(1)1}=-{1\over 2}H^1_{11}y^i}\\
\displaystyle{G^{(i)}_{(1)1}={1\over 2}\gamma^i_{jk}y^j}y^k.
\end{array}\right.
\end{equation}
The equations of the harmonic curves attached to these sprays, with respect to
the semi-Riemannian temporal metric $h$, reduce to
\begin{equation}
h^{11}\left\{{d^2x^i\over dt^2} -H^1_{11}{dx^i\over dt}+\gamma^i_{jk}{dx^j
\over dt}{dx^k\over dt}\right\}=0,
\end{equation}
that is, exactly the equations whose solutions are the well known classical
harmonic maps between the semi-Riemannian manifolds $(R,h)$ and $(M,\varphi)$
\cite{4}. Particularly, if  we regard the temporal manifold $R$ endowed with
the euclidian metric $h=\delta$, we recover the classical equations of geodesics
on the semi-Riemannian manifold $M$. These facts emphasize the naturalness
of our previous definition.

\subsection{Nonlinear connections. Adapted bases.}

\setcounter{equation}{0}
\hspace{5mm} It is well known the importance of the nonlinear connections in the study of
the geometry of a fibre bundle $E$. A nonlinear connection (i. e. a supplementary
horizontal  distribution of the vertical distribution of $E$) offers the
possibility of construction of the {\it vector} or {\it covector adapted bases}.
These allow to write, in a simple form, the geometrical objects or properties of
the total space $E$. In this sense, considering the particular case $E=J^1(R,M)$,
we proved in \cite{16},
\begin{th}
A nonlinear connection $\Gamma$ on the jet fibre bundle of order one $E$
is determined by a pair of local function sets $M^{(i)}_{(1)1}$ and
$N^{(i)}_{(1)j}$ which modify by the transformation laws
\begin{equation}\label{tnc}
\displaystyle{\tilde M^{(j)}_{(1)1}{d\tilde t\over dt}=M^{(k)}_{(1)1}{dt\over
d\tilde t}{\partial\tilde x^j\over\partial x^k}-{\partial\tilde y^j\over
\partial t}},
\end{equation}
\begin{equation}\label{snc}
\displaystyle{\tilde N^{(j)}_{(1)k}{\partial\tilde x^k\over\partial x^i}=
N^{(k)}_{(1)i}{dt\over d\tilde t}{\partial\tilde x^j\over\partial x^k}-
{\partial\tilde y^j\over\partial x^i}}.
\end{equation}
\end{th}
\addtocounter{defin}{1}
{\bf Definition \thedefin} A set of local functions $M^{(i)}_{(1)1}$ (resp.
$N^{(i)}_{(1)j}$) on $J^1(R,M)$, which transform by the rules \ref{tnc} (resp.
\ref{snc})  is called a {\it temporal nonlinear connection} (resp. {\it spatial
nonlinear connection}) on $E=J^1(R,M)$.
\medskip\\
\addtocounter{ex}{1}
{\bf Example \theex} Studying the transformation
rules of the local components
\begin{equation}\label{cnc}
\left\{\begin{array}{l}\medskip
\displaystyle{M^{(i)}_{(1)1}=-H^1_{11}y^i}\\
\displaystyle{N^{(i)}_{(1)j}=\gamma^i_{jk}y^k},
\end{array}\right.
\end{equation}
where $H^1_{11}$ (resp.  $\gamma^i_{jk}$) are the Christoffel symbols of a temporal
(resp. spatial)  semi-Riemannian metric $h$ (resp. $\varphi$), we conclude
that $\Gamma_0=(M^{(i)}_{(1)1},N^{(i)}_{(1)j})$ represents a nonlinear
connection on $E$, which is called the {\it canonical nonlinear connection
attached to the metric pair $(h,\varphi)$.}\medskip

Taking into account the transformation laws \ref{lts}, \ref{lss} and
\ref{tnc}, \ref{snc}, we deduce  without difficulties that the notion
of temporal (resp. spatial) spray is intimately connected to the notion of
temporal (resp. spatial) nonlinear connection.
\begin{th}\label{l1}
i) If $M^{(i)}_{(1)1}$ are the components of a temporal nonlinear connection,
then the components
\begin{equation}
\displaystyle {H^{(i)}_{(1)1}={1\over2}M^{(i)}_{(1)1}}
\end{equation}
represent a temporal spray.

ii) Conversely, if $H^{(i)}_{(1)1}$ are the components of a temporal spray,
then
\begin{equation}
M^{(i)}_{(1)1}=2H^{(i)}_{(1)1}
\end{equation}
are the components of a temporal nonlinear connection.
\end{th}

\begin{th}\label{l2}
i) If $G^{(i)}_{(1)1}$ are the components of a spatial spray, then the
components
\begin{equation}
N^{(i)}_{(1)j}={\partial G^i_{(1)1}\over\partial y^j}
\end{equation}
represent a spatial nonlinear connection.

ii) Conversely, the spatial nonlinear connection $N^{(i)}_{(1)j}$ induces
the spatial spray
\begin{equation}
2G^{(i)}_{(1)1}=N^{(i)}_{(1)j}y^j.
\end{equation}
\end{th}
\addtocounter{rem}{1}
{\bf Remark \therem} The previous theorems allow us to conclude that a
time-dependent spray $(H,G)$  induces naturally a nonlinear connection
$\Gamma$ on $E$, which is called the {\it canonical nonlinear connection
associated to the time-dependent spray $(H,G)$.} We point out that the
canonical nonlinear  connection $\Gamma$ attached  to the time-dependent
spray $(H,G)$ is a natural generalization of the canonical nonlinear
connection $N$ induced by a time-dependent spray $G$ from the classical
rheonomic Lagrangian geometry \cite{11}.\medskip

Let $\Gamma=(M^{(i)}_{(1)1} N^{(i)}_{(1)j})$ be  a nonlinear connection on
the 1-jet fibre bundle $E$.  Let us consider the geometrical objects,
\begin{equation}
\left\{\begin{array}{l}\medskip
\displaystyle{
{\delta\over\delta t}={\partial\over\partial t}-M^{(j)}_{(1)1}{\partial\over
\partial y^j}}\\\medskip
\displaystyle{{\delta\over\delta x^i}={\partial\over\partial x^i}-N^{(j)}_
{(1)i}{\partial\over\partial y^j}}\\
\delta y^i=dy^i+M^{(i)}_{(1)1}dt+N^{(i)}_{(1)j}dx^j.
\end{array}\right.
\end{equation}
One easily deduces that the  set of vector fields
$\displaystyle{\left\{{\delta\over\delta t},{\delta\over\delta
x^i},{\partial\over\partial y^i}\right\}\subset{\cal X}(E)}$ and of covector
fields
$\{dt,dx^i,\delta y^i\}\subset{\cal X}^*(E)$ are  dual bases.\medskip\\
\addtocounter{defin}{1}
{\bf Definition \thedefin} The basis $\displaystyle{\left\{{\delta\over\delta t},
{\delta\over\delta x^i},{\partial\over\partial y^i}\right\}\subset {\cal X}
(E)}$ and its  dual basis $\{dt, dx^i,\delta y^i\}\subset{\cal X}^*(E)$ are
called the {\it adapted bases} on $E$, determined by the nonlinear connection
$\Gamma$.\medskip

The big advantage of the adapted bases is that the transformation laws of its
elements are simple and natural.
\begin{prop}
The transformation laws of the elements of the adapted bases attached to the
nonlinear connection $\Gamma$ are
\begin{equation}\label{vab}
\left\{\begin{array}{l}\medskip
\displaystyle{{\delta\over\delta t}={d\tilde t\over dt}{\delta\over\delta
\tilde t}}\medskip\\
\displaystyle{{\delta\over\delta x^i}={\partial\tilde x^j\over\partial x^i}
{\delta\over\delta\tilde x^j}}\medskip\\
\displaystyle{{\partial\over\partial y^i}={\partial\tilde x^j\over\partial
x^i}{dt\over d\tilde t}{\delta\over\delta\tilde y^j},}
\end{array}\right.
\end{equation}
\begin{equation}\label{cvab}
\left\{\begin{array}{l}\medskip
\displaystyle{dt={dt\over d\tilde t}d\tilde t}\medskip\\
\displaystyle{dx^i={\partial x^i\over\partial\tilde x^j}d\tilde x^j}
\medskip\\
\displaystyle{\delta y^i={\partial x^i\over\partial\tilde x^j}{d\tilde t\over
dt}\delta\tilde y^j.}
\end{array}\right.
\end{equation}
\end{prop}
\addtocounter{rem}{1}
{\bf Remark \therem} The simple transformation rules \ref{vab} and \ref{cvab}
determine us to describe the objects with geometrical and physical meaning
from the subsequent rheonomic Lagrange theory of physical fields, in adapted
components.

\subsection{$\Gamma$-linear connections}

\setcounter{equation}{0}
\hspace{5mm} In order to develope the theory of $\Gamma$-linear connections
on the 1-jet space $E$, we need the following
\begin{prop}
i) The Lie algebra ${\cal X}(E)$ of vector fields decomposes as
$$
{\cal X}(E)={\cal X}({\cal H}_T)\oplus{\cal X}({\cal H}_M)\oplus
{\cal X}({\cal V}),
$$
where
$$
{\cal X}({\cal H}_T)=Span\left\{{\delta\over\delta t^\alpha}\right\},\quad
{\cal X}({\cal H}_M)=Span\left\{{\delta\over\delta x^i}\right\},\quad
{\cal X}({\cal V})=Span\left\{{\partial\over\partial x^i_\alpha}\right\}.
$$

ii) The Lie algebra ${\cal X}^*(E)$ of covector fields decomposes as
$$
{\cal X}^*(E)={\cal X}^*({\cal H}_T)\oplus{\cal X}^*({\cal H}_M)\oplus
{\cal X}^*({\cal V}),
$$
where
$$
{\cal X}^*({\cal H}_T)=Span\{dt^\alpha\},\quad
{\cal X}^*({\cal H}_M)=Span\{dx^i\},\quad
{\cal X}^*({\cal V})=Span\{\delta x^i_\alpha\}.
$$
\end{prop}

Let us consider $h_T$, $h_M$ (horizontal) and $v$ (vertical) as the canonical projections of the above
decompositions.\medskip\\
\addtocounter{defin}{1}
{\bf Definition \thedefin} A linear connection $\nabla:{\cal X}(E)\times{\cal X}(E)\to
{\cal X}(E)$ is called {\it a $\Gamma$-linear connection on $E$} if $\nabla
h_T=0$, $\nabla h_M=0$ and $\nabla v=0$.\medskip

In order to describe in local terms a $\Gamma$-linear connection $\nabla$ on
$E$, we need nine unique local components,
\begin{equation}\label{lglc}
\nabla\Gamma=(\bar G^1_{11},G^k_{i1},G^{(k)(1)}_{(1)(i)1},\bar L^1_{1j},L^k_
{ij},L^{(k)(1)}_{(1)(i)j},\bar C^{1(1)}_{1(j)},C^{k(1)}_{i(j)},C^{(k)(1)(1)}_
{(1)(i)(j)}),
\end{equation}
which are locally defined by the relations
$$
\begin{array}{lll}\medskip
\displaystyle{\nabla_{\delta\over\delta t}{\delta\over\delta t}=\bar G^1_{11}
{\delta\over\delta t},}&\displaystyle{\nabla_{\delta\over\delta t}{\delta
\over\delta x^i}=G^k_{i1}{\delta\over\delta x^k},}&\displaystyle{\nabla_{
\delta\over\delta t}{\partial\over\partial y^i}=G^{(k)(1)}_{(1)(i)1}{\partial
\over\partial y^k},}\\\medskip
\displaystyle{\nabla_{\delta\over\delta x^j}{\delta\over\delta t}=\bar L^1_
{1j}{\delta\over\delta t},}&\displaystyle{\nabla_{\delta\over\delta x^j}
{\delta\over\delta x^i}=L^k_{ij}{\delta\over\delta x^k},}&
\displaystyle{\nabla_{\delta\over\delta x^j}{\partial\over\partial y^i}=
L^{(k)(1)}_{(1)(i)j}{\partial\over\partial y^k},}\\
\displaystyle{\nabla_{\partial\over\partial y^j}{\delta\over\delta t}=\bar
C^{1(1)}_{1(j)}{\delta\over\delta t},}&\displaystyle{\nabla_{\partial\over
\partial y^j}{\delta\over\delta x^i}=C^{k(1)}_{i(j)}{\delta\over\delta x^k},}&
\displaystyle{\nabla_{\partial\over\partial y^j}{\partial\over\partial
y^i}=C^{(k)(1)(1)}_{(1)(i)(j)}{\partial\over\partial y^k}}.
\end{array}
$$
Now, using the transformation laws \ref{vab} of the elements
$\displaystyle{\left\{{\delta\over\delta t},{\delta\over\delta x^i},
{\partial\over\partial y^i}\right\}}$ together with
the properties of the $\Gamma$-linear connection $\nabla$, we obtain by
computations
\begin{th}\label{glc}
i) The coefficients of the $\Gamma$-linear connection $\nabla$ modify
by the rules\\\\
$(h_T)\hspace{6mm}
\left\{\begin{array}{l}\medskip\displaystyle{
\bar G^1_{11}{d\tilde t\over dt}=\tilde{\bar G}^1_{11}\left({d\tilde t\over
dt}\right)^2+{d^2\tilde t\over dt^2}}\\\medskip
\displaystyle{G^k_{i1}=\tilde G^m_{j1}{\partial x^k\over\partial \tilde x^m}
{\partial\tilde x^j\over\partial x^i}{d\tilde t\over dt}}\\
\medskip
\displaystyle{G^{(k)(1)}_{(1)(i)1}=\tilde G^{(m)(1)}_{(1)(j)1}{\partial x^k
\over\partial\tilde x^m}{\partial\tilde x^j\over\partial x^i}{d\tilde t\over
dt}+\delta^k_i\left({d\tilde t\over dt}\right)^2{d^2 t\over d\tilde t^2}},
\end{array}\right.\medskip
$\\
$(h_M)\hspace{5mm}
\left\{\begin{array}{l}\medskip\displaystyle{
\bar L^1_{1j}{\partial x^j\over\partial\tilde x^l}=\tilde{\bar L}^1_{1l}}\\
\medskip
\displaystyle{L^m_{ij}{\partial\tilde x^r\over\partial x^m}=\tilde L^r_{pq}
{\partial\tilde x^p\over\partial x^i}{\partial\tilde x^q\over\partial x^j}+
{\partial^2\tilde x^r\over\partial x^i\partial x^j}}\\\medskip
\displaystyle{L^{(m)(1)}_{(1)(i)j}{\partial\tilde x^r\over\partial x^m}=
\tilde L^{(r)(1)}_{(1)(p)q}{\partial x^p\over\partial\tilde x^i}{\partial
\tilde x^q\over\partial x^j}+{\partial^2\tilde x^r\over\partial x^i\partial
x^j}},
\end{array}\right.\medskip
$\\
$(v)\hspace{8mm}
\left\{\begin{array}{l}\medskip\displaystyle{
\bar C^{1(1)}_{1(i)}=\tilde{\bar C}^{1(1)}_{1(j)}{\partial\tilde x^j\over
\partial x^i}{dt\over d\tilde t}}\\
\medskip
\displaystyle{C^{k(1)}_{i(j)}=\tilde C^{s(1)}_{p(r)}{\partial x^k\over\partial
\tilde x^s}{\partial\tilde x^p\over\partial x^i}{\partial\tilde x^r\over\partial x^j}
{dt\over d\tilde t}}\\
\medskip
\displaystyle{C^{(k)(1)(1)}_{(1)(i)(j)}=\tilde C^{(r)(1)(1)}_{(1)(p)(q)}
{\partial x^k\over\partial\tilde x^r}{\partial\tilde x^p\over\partial x^i}
{\partial\tilde x^q\over\partial x^j}{dt\over d\tilde t}}.
\end{array}\right.\medskip
$

ii) Conversely, to give a $\Gamma$-linear connection $\nabla$ on the 1-jet space
$E$ is equivalent to give a set of nine local coefficients \ref{lglc}
whose local transformations laws are described in i.
\end{th}

The previous theorem allows us to offer an important example of $\Gamma$-linear
connection on  $J^1(R,M)$.\medskip\\
\addtocounter{ex}{1}
{\bf Example \theex} Let $h_{11}$ (resp. $\varphi_{ij}$) be a
semi-Riemannian metric on the temporal (resp. spatial) manifold $R$ (resp.
$M$) and $H^1_{11}$ (resp. $\gamma^k_{ij}$) its Christoffel
symbols. Let us consider $\Gamma_0=(M^{(i)}_{(1)1}, N^{(i)}_{(1)j})$,
where $M^{(i)}_{(1)1}=-H^1_{11}y^i,\;N^{(i)}_{(1)j}=\gamma^i_{jk}y^k$, the
canonical nonlinear connection on $E$ attached to the metric pair $(h_{11},
\varphi_{ij})$. Using the transformation rules of the Christoffel symbols, we
deduce that the following set of local coefficients \cite{15}
\begin{equation}
B\Gamma_0=(\bar G^1_{11},0,G^{(k)(1)}_{(1)(i)1},0,L^k_{ij},L^{(k)(1)}_{(1)(i)j},
0,0,0),
\end{equation}
where $\bar G^1_{11}=H^1_{11},\;G^{(k)(1)}_{(1)(i)1}=-\delta^k_iH^1_{11},\;
L^k_{ij}=\gamma_{ij}^k$ and $L^{(k)(1)}_{(1)(i)j}=\delta^1_1\gamma^k_{ij}$,
is a\linebreak $\Gamma_0$-linear connection. This is called {\it the Berwald
$\Gamma_0$-linear connection of the metric pair $(h_{11},\varphi_{ij})$}.
\medskip

Note that a $\Gamma$-linear connection $\nabla$ on $E$, defined by the local
coefficients \ref{lglc},
induces a natural linear connection on the d-tensors set of the jet fibre
bundle $J^1(R,M)$, in the following fashion.
Starting with $X\in{\cal X}(E)$ a d-vector field and a d-tensor field $D$
locally expressed by
$$\begin{array}{l}\medskip\displaystyle{
X=X^1{\delta\over\delta t}+X^m{\delta\over\delta x^m}
+X^{(m)}_{(1)}{\partial\over\partial y^m},}\\
\displaystyle{
D=D^{1i(j)(1)\ldots}_{1k(1)(l)\ldots}{\delta\over\delta t}\otimes{\delta\over
\delta x^i}\otimes{\partial\over\partial y^j}\otimes dt\otimes dx^k\otimes
\delta y^l\ldots,}
\end{array}
$$
we introduce the covariant derivative
$$
\begin{array}{l}\medskip
\nabla_XD=X^1\nabla_{\delta\over\delta t}D+X^p\nabla_{\delta\over\delta x^p}D
+X^{(p)}_{(1)}\nabla_{\partial\over\partial y^p}D=\left\{X^1D^{1i(j)(1)\ldots}
_{1k(1)(l)\ldots /1}+X^p\right.\\\left.
D^{1i(j)(1)\ldots}_{1k(1)(l)\ldots\vert p}+X^{(p)}_{(1)}D^{1i(j)(1)\ldots}
_{1k(1)(l)\ldots}\vert_{(p)}^{(1)}\right\}\displaystyle{{\delta\over\delta
t}\otimes{\delta\over\delta x^i}\otimes{\partial\over\partial y^j}\otimes
dt\otimes dx^k\otimes\delta y^l\ldots,}
\end{array}
$$
where\\\\
$
(h_T)\hspace{6mm}\left\{\begin{array}{l}\medskip\displaystyle{
D^{1i(j)(1)\ldots}_{1k(1)(l)\ldots /1}={\delta D^{1i(j)(1)\ldots}_{1k(1)(l)
\ldots}\over\delta t}+D^{1i(j)(1)\ldots}_{1k(1)(l)\ldots}\bar G^1_{11}+}\\
\medskip
+D^{1m(j)(1)\ldots}_{1k(1)(l)\ldots}G^i_{m1}+D^{1i(m)(1)\ldots}_{1k(1)(l)
\ldots}G^{(j)(1)}_{(1)(m)1}+\ldots-\\
-D^{1i(j)(1)\ldots}_{1k(1)(l)\ldots}\bar G^1_{11}-D^{1i(j)(1)\ldots}_{1m(1)(l)
\ldots} G^m_{k1}-D^{1i(j)(1)\ldots}_{1k(1)(m)\ldots} G^{(m)(1)}_{(1)(l)1}-
\ldots,
\end{array}\right.\medskip
$
$
(h_M)\hspace{5mm}\left\{\begin{array}{l}\medskip\displaystyle{
D^{1i(j)(1)\ldots}_{1k(1)(l)\ldots\vert p}={\delta D^{1i(j)(1)\ldots}_{1k(1)
(l)\ldots}\over\delta x^p}+D^{1i(j)(1)\ldots}_{1k(1)(l)\ldots}\bar L^1_{1p}+}
\\\medskip
+D^{1m(j)(1)\ldots}_{1k(1)(l)\ldots}L^i_{mp}+D^{1i(m)(1)\ldots}_{1k(1)(l)
\ldots}L^{(j)(1)}_{(1)(m)p}+\ldots-\\
-D^{1i(j)(1)\ldots}_{1k(1)(l)\ldots}\bar L^1_{1p}-D^{1i(j)(1)\ldots}_{1m(1)
(l)\ldots} L^m_{kp}-D^{1i(j)(1)\ldots}_{1k(1)(m)\ldots} L^{(m)(1)}_{(1)(l)p}-
\ldots,
\end{array}\right.\medskip
$
$
(v)\hspace{8mm}\left\{\begin{array}{l}\medskip\displaystyle{
D^{1i(j)(1)\ldots}_{1k(1)(l)\ldots}\vert_{(p)}^{(1)}={\partial D^{1i(j)(1)
\ldots}_{1k(1)(l)\ldots}\over\partial y^p}+D^{1i(j)(1)\ldots}_{1k(1)(l)\ldots}
\bar C^{1(1)}_{1(p)}+}\\\medskip
+D^{1m(j)(1)\ldots}_{1k(1)(l)\ldots}C^{i(1)}_{m(p)}+D^{1i(m)(1)\ldots}_{1k(1)
(l)\ldots}C^{(j)(1)(1)}_{(1)(m)(p)}+\ldots-\\
-D^{1i(j)(1)\ldots}_{1k(1)(l)\ldots}\bar C^{1(1)}_{1(p)}-D^{1i(j)(1)\ldots}_
{1m(1)(l)\ldots} C^{m(1)}_{k(p)}-D^{1i(j)(1)\ldots}_{1k(1)(m)\ldots} C^{(m)
(1)(1)}_{(1)(l)(p)}-\ldots.
\end{array}\right.\medskip
$

The local operators "$ _{/1}$", "$_{\vert p}$" and "$\vert^{(1)}_{(p)}$" are
called the {\it $H_R$-horizontal covariant derivative, $h_M$-horizontal covariant
derivative} and {\it $v$-vertical covariant derivative} of the $\Gamma$-linear
connection $\nabla$.

The study of the torsion {\bf T} and curvature {\bf R} d-tensors of an arbitrary
$\Gamma$-linear connection $\nabla$ was made  in \cite{18}. In this
context, we proved that the torsion\linebreak d-tensor is
determined by twelve effective local torsion d-tensors, while the curvature
d-tensor of $\nabla$ is determined by eighteen local d-tensors.

\subsection{$h$-Normal $\Gamma$-linear connections}

\setcounter{equation}{0}
\hspace{5mm} Let $h_{11}$ be a fixed pseudo-Riemannian metric on the temporal
manifold $R$, $H^1_{11}$ its Christoffel symbols and $J=J^{(i)}_{(1)1j}
{\partial\over\partial y^i}\otimes dt\otimes dx^j$, where $J^{(i)}_{(1)1j}=
h_{11}\delta^i_j$, the {\it normalization d-tensor} \cite{16} attached to the metric
$h_{11}$. In order  to reduce the big number of torsion and curvature d-tensors which characterize
a general $\Gamma$-linear connection on $E$, we consider the
following\medskip\\
\addtocounter{defin}{1}
{\bf Definition \thedefin} A $\Gamma$-linear connection $\nabla$ on $E$,
defined by the local coefficients
$$
\nabla\Gamma=(\bar G^1_{11},G^k_{i1},G^{(k)(1)}_{(1)(i)1},\bar L^1_{1j},L^k_
{ij},L^{(k)(1)}_{(1)(i)j},\bar C^{1(1)}_{1(j)},C^{k(1)}_{i(j)},C^{(k)(1)(1)}_
{(1)(i)(j)}),
$$
that verify the relations $\bar G^1_{11}=H^1_{11},\;\bar L^1_{1j}=0,\;\bar
C^{1(1)}_{1(j)}=0$ and $\nabla J=0$, is called a {\it $h$-normal $\Gamma$-linear
connection}.\medskip\\
\addtocounter{rem}{1}
{\bf Remark \therem} Taking into account the local covariant $h_R$-horizontal
"$_{/1}$", $h_M$-horizontal "$_{\vert k}$" and $v$-vertical "$\vert^{(1)}_{(k)}$"
covariant derivatives induced by $\nabla$, the condition $\nabla J=0$ is
equivalent to
\begin{equation}
J^{(i)}_{(1)1j/1}=0,\quad
J^{(i)}_{(1)1j\vert k}=0,\quad
J^{(i)}_{(1)1j}\vert^{(1)}_{(k)}=0.
\end{equation}

In this context, we can prove the following

\begin{th}
The coefficients of a $h$-normal $\Gamma$-linear
connection $\nabla$ verify the identities
\begin{equation}
\begin{array}{lll}\medskip
\bar G^1_{11}=H^1_{11},&\bar L^1_{1j}=0,&
\bar C^{1(1)}_{1(j)}=0,\\\medskip
G^{(k)(1)}_{(1)(i)1}=G^k_{i1}-\delta^k_iH^1_{11},&L^{(k)(1)}_{(1)(i)j}=L^k_
{ij},&C^{(k)(1)(1)}_{(1)(i)(j)}=C^{k(1)}_{i(j)}.
\end{array}
\end{equation}
\end{th}
{\bf Proof.} The first three relations come from the definiton of a $h$-normal
$\Gamma$-linear connection.

The condition $\nabla J=0$ implies locally that
\begin{equation}
\left\{\begin{array}{l}\medskip
\displaystyle{h_{11}G^{(i)(1)}_{(1)(j)1}=h_{11}G^i_{j1}+\delta^i_j\left[-
{\partial h_{11}\over\partial t}+H_{111}\right]}\\\medskip
h_{11}L^{(i)(1)}_{(1)(j)}=h_{11}L^i_{jk}\\
h_{11}C^{(i)(1)(1)}_{(1)(j)(k)}=h_{11}C^{i(1)}_{j(k)},
\end{array}\right.
\end{equation}
where $H_{111}=H^1_{11}h_{11}$ represent the Christoffel symbols of the first
kind attached to the semi-Riemannian metric $h_{11}$. Contracting the above
relations by $h^{11}$, one obtains the last three identities of the theorem.
\rule{5pt}{5pt}\medskip\\
\addtocounter{rem}{1}
{\bf Remarks \therem} i) The preceding theorem implies that a $h$-normal
$\Gamma$-linear on $E$ is determined just by four effective coefficients
$$
\nabla\Gamma=(H^1_{11},G^k_{i1},L^k_{ij},C^{k(1)}_{i(j)}).
$$

ii) Considering the particular case  of the temporal metric $h=\delta$, we
remark that a $\delta$-normal $\Gamma$-linear connection on $J^1(R,M)$ is a
natural generalization of the notion of $N$-linear connection used in the
\cite{11}.\medskip\\
\addtocounter{ex}{1}
{\bf Example \theex} Using the previous theorem, we deduce that the canonical
Berwald $\Gamma_0$-linear connection associated to the metric pair $(h_{11},
\varphi_{ij})$ is a $h$-normal $\Gamma_0$-linear connection, defined by the
local coefficients $B\Gamma_0=(H^1_{11},0,\gamma^k_{ij},0)$.\medskip

\subsection{d-Torsions and d-Curvatures}

\setcounter{equation}{0}
\hspace{5mm} The study of the torsion {\bf T} and curvature {\bf R} d-tensors
of an arbitrary $h$-normal $\Gamma$-linear connection $\nabla$ was made in
\cite{15}. We proved there that  the adapted components $\bar T^1_{11},\;
\bar T^1_{1j}$, $\bar P^{1(1)}_{1(j)}$ and $R^{(m)}_{(1)11}$ of the torsion
d-tensor {\bf T} of $\nabla$ vanish. Consequently, we obtain the following
\cite{15}
\begin{th}
The torsion d-tensor {\bf T} of the $h$-normal $\Gamma$-linearconnection
$\nabla$ is determined by eight local d-tensors
\begin{equation}
\begin{tabular}{|c|c|c|c|}
\hline
&$h_T$&$h_M$&$v$\\
\hline
$h_Th_T$&0&0&$0$\\
\hline
$h_Mh_T$&0&$T^m_{1j}$&$R^{(m)}_{(1)1j}$\\
\hline
$h_Mh_M$&0&$T^m_{ij}$&$R^{(m)}_{(1)ij}$\\
\hline
$vh_T$&0&0&$P^{(m)\;\;(1)}_{(1)1(j)}$\\
\hline
$vh_M$&0&$P^{m(1)}_{i(j)}$&$P^{(m)\;(1)}_{(1)i(j)}$\\
\hline
$vv$&0&0&$S^{(m)(1)(1)}_{(1)(i)(j)}$\\
\hline
\end{tabular}
\end{equation}
where
$
\displaystyle{P^{(m)\;\;(1)}_{(1)1(j)}={\partial M^{(m)}_{(1)1}\over\partial
y^j}-G^m_{j1}+\delta^m_jH^1_{11},}\quad
$
$\displaystyle{P^{(m)\;\;(1)}_{(1)i(j)}={\partial N^{(m)}_{(1)i}\over\partial
y^j}-L^m_{ji},}
$\medskip\linebreak
$\displaystyle{
R^{(m)}_{(1)1j}={\delta M^{(m)}_{(1)1}\over\delta x^j}-{\delta N^{(m)}_{(1)j}
\over\delta t},}
$
$\displaystyle{
R^{(m)}_{(1)ij}={\delta N^{(m)}_{(1)i}\over\delta x^j}-{\delta N^{(m)}_{(1)j}
\over\delta x^i},}
$
$
S^{(m)(1)(1)}_{(1)(i)(j)}=C^{m(1)}_{i(j)}-C^{m(1)}_{j(i)}
$,\medskip\linebreak
$
T^m_{1j}=-G^m_{j1}$,\quad
$T^m_{ij}=L^m_{ij}-L^m_{ji}$,\quad
$P^{m(1)}_{i(j)}=C^{m(1)}_{i(j)}$.
\end{th}
\addtocounter{rem}{1}
{\bf Remark \therem} For the Berwald $\Gamma_0$-linear connection associated to the
metrics $h_{11}$ and $\varphi_{ij}$, all torsion d-tensors vanish,
except $R^{(m)}_{(\mu)ij}=r^m_{ijl}x^l_\mu,$ where (resp. $r^m_{ijl}$) are
the curvature tensors of the metric $\varphi_{ij}$.\medskip

In the same context, following the paper \cite{15},  we deduce that the
number of the effective adapted components of the curvature d-tensor {\bf R}
of an $h$-normal $\Gamma$-linear connection $\nabla$ is five.
\begin{th}
The curvature d-tensor {\bf  R} of $\nabla$ is determined by the following
effective local  d-curvatures
\begin{equation}
\begin{tabular}{|c|c|c|c|}
\hline
&$h_T$&$h_M$&$v$\\
\hline
$h_Th_T$&$0$&$0$&$0$\\
\hline
$h_Mh_T$&0&$R^l_{i1k}$&$R^{(l)(1)}_{(1)(i)1k}=R^l_{i1k}$\\
\hline
$h_Mh_M$&0&$R^l_{ijk}$&$R^{(l)(1)}_{(1)(i)jk}=R^l_{ijk}$\\
\hline
$vh_T$&0&$P^{l\;\;(1)}_{i1(k)}$&$P^{(l)(1)\;\;(1)}_{(1)(i)1(k)}=P^{l\;\;(1)}_
{i1(k)}$\\
\hline
$vh_M$&0&$P^{l\;(1)}_{ij(k)}$&$P^{(l)(1)\;(1)}_{(1)(i)j(k)}=P^{l\;(1)}_{ij(k)}$\\
\hline
$vv$&0&$S^{l(1)(1)}_{i(j)(k)}$&$S^{(l)(1)(1)(1)}_{(1)(i)(j)(k)}=S^{l(1)(1)}_
{i(j)(k)}$\\
\hline
\end{tabular}
\end{equation}
where\medskip

$
\displaystyle{R^l_{i1k}={\delta G^l_{i1}\over\delta x^k}-{\delta L^l_{ik}
\over\delta t}+G^m_{i1}L^l_{mk}-L^m_{ik}G^l_{m1}+C^{l(1)}_{i(m)}R^{(m)}_{(1)1k},}
$\medskip

$
\displaystyle{R^l_{ijk}={\delta L^l_{ij}\over\delta x^k}-
{\delta L^l_{ik}\over\delta x^j}+L^m_{ij}L^l_{mk}-L^m_{ik}L^l_{mj}+
C^{l(1)}_{i(m)}R^{(m)}_{(1)jk},}
$\medskip

$
\displaystyle{P^{l\;\;(1)}_{i1(k)}={\partial G^l_{i1}\over\partial y^k}-C^
{l(1)}_{i(k)/1}+C^{l(1)}_{i(m)}P^{(m)\;\;(1)}_{(1)1(k)},}
$\medskip

$
\displaystyle{P^{l\;(1)}_{ij(k)}={\partial L^l_{ij}\over\partial y^k}-C^
{l(1)}_{i(k)\vert j}+C^{l(1)}_{i(m)}P^{(m)\;(1)}_{(1)j(k)},}
$\medskip

$
\displaystyle{S^{l(1)(1)}_{i(j)(k)}={\partial C^{l(1)}_{i(j)}\over\partial
y^k}-{\partial C^{l(1)}_{i(k)}\over\partial y^j}+C^{m(1)}_{i(j)}C^{l(1)}_
{m(k)}-C^{m(1)}_{i(k)}C^{l(1)}_{m(j)}.}
$
\end{th}
\addtocounter{rem}{1}
{\bf Remark \therem} In the case of the Berwald $\Gamma_0$-linear connection
associated to the metric pair $(h_{11},\varphi_{ij})$, all curvature
d-tensors vanish, except $R^l_{ijk}=r^l_{ijk}$, where $r^l_{ijk}$ are the
curvature tensors of the metric $\varphi_{ij}$.\pagebreak

\section{Relativistic rheonomic Lagrange geometry}

\subsection{Some aspects of classical rheonomic Lagrange geometry}
\setcounter{equation}{0}
\hspace{5mm} A lot of geometrical models in Mechanics, Physics or Biology are
based on the notion of ordinary Lagrangian. Thus, the concept of Lagrange space which
generalizes that of Finsler space was introduced. In order to
geometrize the fundamental concept in mechanics, that of Lagrangian, we recall
that a Lagrange  space $L^n=(M,L(x,y))$ is defined as a pair which consists
of a real, smooth, $n$-dimensional manifold $M$ and a regular Lagrangian
$L:TM\to R$, not necessarily homogenous with respect to the direction $(y^i)_
{i=\overline{1,n}}$. The differential geometry of Lagrange spaces is now
considerably developped and used in various fields to study natural process
where the dependence on position, velocity  or momentum is involved \cite{11}.
Also, the geometry of Lagrange spaces gives a model for both the gravitational
and electromagnetic field, in a very natural blending of the geometrical
structure of the space with the characteristic properties of these physical
fields.

At the same time, there are many problems in Physics and  Variational calculus
in which time dependent Lagrangians (i. e., a  smooth real function on
$R\times TM$) are involved. A geometrization of a such time dependent
Lagrangian is sketched in \cite{11}. This is called the {\it "Rheonomic
Lagrange Geometry"}. On the one hand, it is remarkable that this geometrical
model is the house of the development of the {\it classical rheonomic
Lagrangian mechanics}. On the other hand, from our point of view, this time
dependent Lagrangian geometrization has an important inconvenience that
we will describe.

In the context exposed in the book \cite{11}, the
energy action functional ${\cal E}$, attached to a given time dependent
Lagrangian,
$$
L:R\times TM\to R,\quad (t,x^i,y^i)\to L(t,x^i,y^i),
$$
not necessarily homogenous with respect to the direction
$(y^i)_{i=\overline{1,n}}$, is of the form
\begin{equation}
{\cal E}(c)=\int^b_aL(t,x^i(t),\dot x^i(t))\;dt,
\end{equation}
where $[a,b]\subset R$, and $c:[a,b]\to M$ is a smooth curve, locally expressed
by $t\to(x^i(t))$, and having the velocity $\dot x=(\dot x^i(t))$.
It is obvious that the non-homogeneity of the Lagrangian $L$, regarded as a smooth
function on the product manifold $R\times TM$, implies that the energy
action functional ${\cal E}$ is {\it dependent of the parametrizations of
every curve $c$}. In order to
remove this difficulty,  the authors regard the space $R\times TM$  like a
fibre bundle over $M$. In this context,  the geometrical invariance group of
$R\times TM$ is  given by \ref{G_2}. In other words, to remove  the parametrization
dependence of ${\cal E}$, they ignore the temporal repametrizations on $R\times TM$.
Naturally, in these conditions, their energy functional becomes a well defined
one, but their approach stands out by the  {\it  "absolute"} character of the
time $t$.

In our geometrical approach, we try to remove this inconvenience. For that
reason we regard the space $R\times TM\equiv J^1(R,M)$ like a fibre bundle
over $R\times M$. The gauge group of this bundle of configurations is given by
\ref{G_1}. Consequently, our gauge group does not ignore the temporal
reparametrizations, hence, it stands  out by the {\it relativistic} character
of the time $t$. In these conditions, using a given semi-Riemannian metric
$h_{11}(t)$ on $R$, we construct the more general and natural energy action
functional, setting
\begin{equation}
{\cal E}(c)=\int^b_aL(t,x^i(t),\dot x^i(t))\sqrt{\vert h_{11}\vert}\;dt.
\end{equation}
Obviously, ${\cal E}$ is  well defined and is {\it independent of the curve
parametrizations}.

In conclusion, we consider that the difficulty arised in the classical rheonomic
geometry, comes from a puzzling utilization of the notion of Lagrangian.
From this point  of view, we point out that, in our geometrical development,
we use the distinct notions:

i) {\it time dependent Lagrangian function} $-$ A smooth function on $J^1(R,M)$;

ii) {\it time dependent Lagrangian} (Olver's terminology) $-$ A local function ${\cal L}$
on $J^1(R,M)$, which transforms by the rule $\tilde{\cal L}={\cal L}\vert dt/
d\tilde t\vert$. If $L$ is a Lagrangian function on 1-jet fibre bundle,  then
${\cal L}=L\sqrt{\vert  h_{11}\vert}$ represents a Lagrangian on $J^1(R,M)$.

Finally, we point out that the geometrization attached to a time-dependent
Lagrangian function that we will construct, can be called
{\it "Relativistic Rheonomic Lagrange Geometry"}. From our  point of view, this
geometry becomes a natural instrument in the development of the {\it relativistic
rheonomic Lagrangian mechanics}.

\subsection{Relativistic rheonomic Lagrange spaces}

\setcounter{equation}{0}
\hspace{5mm} In order to develope our time-dependent Lagrange geometry, we
start the study considering $L:E\to R$
a smooth Lagrangian function on $E=J^1(R,M)$, which is locally expressed by
$E\ni(t,x^i,y^i)\to L(t,x^i,y^i)\in R$. The {\it vertical fundamental metrical
d-tensor} of $L$ is defined by
\begin{equation}
G^{(1)(1)}_{(i)(j)}={1\over 2}{\partial^2L\over\partial y^i\partial y^j}.
\end{equation}

Let $h=(h_{11})$ be a semi-Riemannian metric on the temporal manifold $R$.\medskip\\
\addtocounter{defin}{1}
{\bf Definition \thedefin} A Lagrangian function $L:E\to R$ whose vertical
fundamental metrical d-tensor is of the form
\begin{equation}
G^{(1)(1)}_{(i)(j)}(t,x^k,y^k)=h^{11}(t)g_{ij}(t,x^k,y^k),
\end{equation}
where $g_{ij}(t,x^k,y^k)$ is a  d-tensor on $E$,  symmetric,
of rank $n$  and  having a  constant signature on $E$, is called a {\it
Kronecker $h$-regular Lagrangian function, with respect to the temporal
semi-Riemannian metric $h=(h_{11})$}.

In this context, we can introduce the  following\medskip\\
\addtocounter{defin}{1}
{\bf  Definition \thedefin} A pair $RL^n=(J^1(R,M),L)$, where $n=\dim M$,
which consists of the  1-jet  fibre bundle and a
Kronecker $h$-regular Lagrangian function \linebreak
$L:J^1(T,M)\to R$ is called a
{\it relativistic  rheonomic Lagrange space}.\medskip\\
\addtocounter{rem}{1}
{\bf Remark \therem} In our geometrization of the time-dependent Lagrangian
function $L$ that we will construct, all entities with geometrical or
physical meaning will be directly arised from the vertical fundamental metrical
d-tensor $G^{(1)(1)}_{(i)(j)}$. This fact points out the {\it metrical
character} (see \cite{5}) and the naturalness of the subsequent relativistic rheonomic
Lagrangian geometry.\medskip\\
\addtocounter{ex}{1}
{\bf Examples \theex} i) Suppose that the spatial manifold $M$ is also endowed with
a semi-Riemannian metric $g=(g_{ij}(x))$. Then, the time dependent Lagrangian function
$L_1:J^1(R,M)\to R$ defined by
\begin{equation}
L_1=h^{11}(t)g_{ij}(x)y^i y^j
\end{equation}
is  a Kronecker $h$-regular time dependent Lagrangian function. Consequently, the pair
$RL^n=(J^1(R,M),L_1)$ is
a relativistic rheonomic Lagrange space. We underline that the Lagrangian
${\cal L}_1=L_1\sqrt{\vert h_{11}\vert}$ is exactly the energy Lagrangian
whose extremals are the harmonic maps between the semi-Riemannian manifolds
$(R,h)$ and $(M,g)$. At the same time, this Lagrangian is a basic object
in the physical theory of bosonic strings.

ii) In above notations, taking $U^{(1)}_{(i)}(t,x)$ as a d-tensor field on
$E$ and $F:R\times M\to R$ a smooth map, the more general Lagrangian function
$L_2:E\to R$ defined by
\begin{equation}
L_2=h^{11}(t)g_{ij}(x)y^i y^j+U^{(1)}_{(i)}(t,x)y^i+F(t,x)
\end{equation}
is also a Kronecker $h$-regular Lagrangian. The relativistic rheonomic
Lagrange space $RL^n=(J^1(R,M),L_2)$ is called the {\it autonomous relativistic
rheonomic Lagrange space  of electrodynamics} because, in the particular case
$h_{11}=1$, we recover the classical Lagrangian space of electrodynamics
\cite{11} which governs the movement law  of a particle placed concomitantly into
a gravitational field and an electromagnetic one. From a physical point of view,
the semi-Riemannian metric $h_{11}(t)$ (resp. $g_{ij}(x)$) represents the {\it gravitational
potentials} of the space $R$ (resp. $M$), the d-tensor $U^{(1)}_{(i)}(t,x)$
stands for the {\it electromagnetic potentials} and $F$ is a function which is
called {\it potential function}. The  non-dynamical character of spatial gravitational
potentials $g_{ij}(x)$ motivates us to use the term of {\it "autonomous"}.

iii) More general, if  we consider $g_{ij}(t,x)$ a d-tensor field  on $E$,
symmetric,  of rank $n$ and  having a constant signature on $E$, we can
define the Kronecker $h$-regular Lagrangian function $L_3:E\to R$,
setting
\begin{equation}
L_3=h^{11}(t)g_{ij}(t,x)y^i y^j+U^{(1)}_{(i)}(t,x)y^i+F(t,x).
\end{equation}
The pair $RL^n=(J^1(R,M),L_3)$ is a relativistic rheonomic Lagrange  space which
is called the
{\it non-autonomous relativistic rheonomic Lagrange space of electrodynamics}.
Physically, we
remark that the gravitational potentials $g_{ij}(t,x)$ of the spatial manifold
$M$ are dependent of the temporal coordinate $t$, emphasizing their dynamic
character.

\subsection{Canonical nonlinear connection}

\setcounter{equation}{0}
\hspace{5mm} Let us consider $h=(h_{11})$ a fixed semi-Riemannian metric
on $R$ and a rheonomic Lagrange space $RL^n=(J^1(R,M),L)$, where $L$ is  a
Kronecker $h$-regular Lagrangian function. Let  $[a,b]\subset R$ be a compact
interval in the temporal manifold $R$. In this context, we can define the {
\it energy action functional} of $RL^n$, setting
$$
{\cal E}:C^\infty(R,M)\to R,\quad
{\cal E}(c)=\int_a^bL(t,x^i,y^i)\sqrt{\vert h\vert}dt,
$$
where the smooth curve $c$ is locally expressed by $(t)\to(x^i(t))$
and $\displaystyle{y^i={dx^i\over dt}}$.

The extremals of the energy functional ${\cal E}$ verifies the Euler-Lagrange
equations
\begin{equation}\label{el}
\hspace*{5mm}2G^{(1)(1)}_{(i)(j)}{d^2x^j\over dt^2}+{\partial^2L\over\partial x^j
\partial y^i}{dx^j\over dt}-{\partial L\over\partial x^i}+{\partial^2L\over
\partial t\partial y^i}+{\partial L\over\partial y^i}H^1_{11}=0,\quad
\forall\;i=\overline{1,n},
\end{equation}
where $H^1_{11}$ are the Christoffel symbols of the semi-Riemannian metric
$h_{11}$.

Taking into account the Kronecker  $h$-regularity of the Lagrangian function
$L$, it is  possible  to rearrange the Euler-Lagrange equations \ref{el} of
the Lagrangian\linebreak ${\cal L}=L\sqrt{\vert h\vert}$, in the Poisson
form \cite{16}
\begin{equation}\label{harm}
\Delta_hx^k+2{\cal G}^k(t,x^m,y^m)=0,\quad\forall\;k=\overline{1,n},
\end{equation}
where
\begin{equation}
\begin{array}{l}\medskip
\displaystyle{\Delta_hx^k=h^{11}\left\{{d^2x^k\over dt^2}-H^1_{11}
{dx^k\over dt}\right\},\;y^m={dx^m\over dt}},\quad\\
\displaystyle{2{\cal G}^k={g^{ki}\over 2}\left\{{\partial^2L\over\partial x^j\partial
y^i}y^j-{\partial L\over\partial x^i}+{\partial^2L\over\partial t\partial y^i}
+{\partial L\over\partial y^i}H^1_{11}+2g_{ij}h^{11}H^1_{11}y^j\right\}.}
\end{array}
\end{equation}
\begin{th}
Denoting $G^{(r)}_{(1)1}=h_{11}{\cal G}^r$, the geometrical object
$G=(G^{(r)}_{(1)1})$ is a spatial spray on the 1-jet space $E$.
\end{th}
{\bf Proof.} By a  direct calculation, we deduce that the local geometrical entities of
the 1-jet space $J^1(R,M)$
\begin{equation}
\begin{array}{l}\medskip
\displaystyle{2{\cal S}^k={g^{ki}\over 2}\left\{{\partial^2L\over\partial x^j
\partial y^i}y^j-{\partial L\over\partial x^i}\right\}}\\\medskip
\displaystyle{2{\cal H}^k={g^{ki}\over 2}\left\{{\partial^2L\over\partial t
\partial y^i}+{\partial L\over\partial y^i}H^1_{11}\right\}}\\
2{\cal J}^k=h^{11}H_{11}^1y^j
\end{array}
\end{equation}
verify the following transformation rules
\begin{equation}
\begin{array}{l}\medskip
\displaystyle{2{\cal S}^p=2\tilde{\cal S}^r{\partial x^p\over\partial\tilde x^r}+h^{11}
{\partial x^p\over\partial\tilde x^l}{d\tilde t\over dt}
{\partial\tilde x^l_\gamma\over\partial x^j}y^j}\\\medskip
\displaystyle{2{\cal H}^p=2\tilde{\cal H}^r{\partial x^p\over\partial\tilde x^r}+h^{11}
{\partial x^p\over\partial\tilde x^l}{d\tilde t\over dt}
{\partial\tilde y^l\over\partial t}}\\
\displaystyle{2{\cal J}^p=2\tilde{\cal J}^r{\partial x^p\over\partial\tilde x^r}-h^{11}
{\partial x^p\over\partial\tilde x^l}{d\tilde t\over dt}
{\partial\tilde y^l\over\partial t}}.
\end{array}
\end{equation}
Consequently, the local entities $2{\cal G}^p=2{\cal S}^p+2{\cal H}^p+2{\cal J}^p$
modify by the transformation laws
\begin{equation}\label{hss}
2\tilde{\cal G}^r=2{\cal G}^p{\partial\tilde x^r\over\partial x^p}-h^{11}
{\partial x^p\over\partial\tilde x^j}{\partial\tilde x^r_\mu\over\partial x^p}
\tilde y^j.
\end{equation}
Hence, multiplying the relation \ref{hss} by $h_{11}$ and regarding the equations
\ref{lss}, we obtain what we were looking for. \rule{5pt}{5pt}\medskip

Taking into account the harmonic curve equations \ref{hc} of a time-dependent
spray on $E$, we can give the  following natural geometrical interpretation
of the Euler-Lagrange equations \ref{harm} attached to the Lagrangian
${\cal L}$:
\begin{th}
The extremals of the energy functional attached to a Kronecker\linebreak $h$-regular
Lagrangian function  $L$ on $J^1(R,M)$ are harmonic curves of the time-dependent
spray $(H,G)$, with respect to the  semi-Riemannian metric $h$, defined by the temporal
components
\begin{equation}
H^{(i)}_{(1)1}=-{1\over 2}H^1_{11}(t)y^i
\end{equation}
and the local spatial components
\begin{equation}
\hspace*{5mm}G^{(i)}_{(1)1}={h_{11}g^{ik}\over 4}\left[{\partial^2L\over\partial x^j
\partial y^k}y^j-{\partial L\over\partial x^k}+{\partial^2L\over\partial
t\partial y^k}+{\partial L\over\partial x^k}H^1_{11}+2h^{11}H^1_{11}g_{kl}
y^l\right].
\end{equation}
\end{th}
\addtocounter{defin}{1}
{\bf Definition \thedefin} The time-dependent spray $(H,G)$ constructed from
the previous theorem is called the {\it canonical time-dependent  spray
attached to the relativistic rheonomic Lagrange space  $RL^n$.}\medskip\\
\addtocounter{rem}{1}
{\bf Remark \therem} In the particular case of an autonomous electrodynamics
relativistic rheonomic Lagrange space (i. e., $g_{ij}(t,x^k,y^k)=g_{ij}(x^k)$),
the canonical spatial spray $G$ is given  by the components
\begin{equation}
G^{(i)}_{(1)1}={1\over 2}\gamma^i_{jk}y^jy^k+
{h_{11}g^{li}\over 4}\left[U^{(1)}_{(l)j}y^j+{\partial U^{(1)}_{(l)}\over
\partial t}+U^{(1)}_{(l)}H^1_{11}-{\partial F\over\partial x^l}\right],
\end{equation}
where
$
\displaystyle{U^{(1)}_{(i)j}={\partial U^{(1)}_{(i)}\over\partial x^j}-
{\partial U^{(1)}_{(j)}\over\partial x^i}}.
$\medskip

In the sequel, using the theorems \ref{l1} and \ref{l2}, we obtain the following
\begin{th}
The pair of local functions $\Gamma=(M^{(i)}_{(1)1},N^{(i)}_{(1)j})$, which
consists of the temporal components
\begin{equation}
M^{(i)}_{(1)1}=2H^{(i)}_{(1)1}=-H^1_{11}y^i,
\end{equation}
and the spatial components
\begin{equation}
N^{(i)}_{(1)j}={\partial G^i_{(1)1}\over\partial y^j},
\end{equation}
where $H^{(i)}_{(1)1}$ and $G^{(i)}_{(1)1}$ are the components  of the canonical
time-dependent spray of $RL^n$,  represents a nonlinear connection on $J^1(R,M)$.
\end{th}
\addtocounter{defin}{1}
{\bf Definition \thedefin} The nonlinear connection $\Gamma=(M^{(i)}_{(1)1},
N^{(i)}_{(1)j})$ from the preceding theorem is called the {\it canonical
nonlinear connection of the relativistic rheonomic Lagrange space $RL^n$.}
\medskip\\
\addtocounter{rem}{1}
{\bf Remark \therem} i) In the case of  an autonomous electrodynamics relativistic
rheonomic Lagrange space (i. e., $g_{ij}(t,x^k,y^k)=g_{ij}(x^k)$), the canonical
nonlinear connection becomes $\Gamma=(M^{(i)}_{(1)1},N^{(i)}_{(1)j})$,
where
\begin{equation}
M^{(i)}_{(1)1}=-H^1_{11}y^i,\quad
N^{(i)}_{(1)j}=\gamma^i_{jk}y^k+{h_{11}g^{ik}\over 4}U^{(1)}_{(k)j}.
\end{equation}

\subsection{Cartan canonical metrical connection}

\setcounter{equation}{0}
\hspace{5mm} The main theorem of this paper  is the theorem of existence of the {\it
Cartan canonical $h$-normal linear connection} $C\Gamma$ which  allow the
subsequent development of the  {\it relativistic rheonomic Lagrangian geometry
of physical fields}, which will be exposed in the next Sections.

\begin{th}
(of existence and uniqueness of Cartan canonical connection)\\
On the relativistic rheonomic Lagrange space $RL^n=(J^1(R,M),L)$ endowed with
its canonical nonlinear connection $\Gamma$ there is a unique $h$-normal
$\Gamma$-linear connection
$$
C\Gamma=(H^1_{11},G^k_{j1},L^i_{jk},C^{i(1)}_{j(k)})
$$
having the metrical properties\medskip

i) $g_{ij\vert k}=0,\quad g_{ij}\vert^{(1)}_{(k)}=0$,\medskip

ii) $\displaystyle{ G^k_{j1}={g^{ki}\over 2}{\delta g_{ij}\over\delta t},
\quad L^k_{ij}=L^k_{ji},\quad C^{i(1)}_{j(k)}=C^{i(1)}_{k(j)}}$.
\end{th}

{\bf Proof.} Let $C\Gamma=(\bar G^1_{11},G^k_{j1},L^i_{jk},C^{i(1)}_{j(k)})$
be a h-normal $\Gamma$-linear connection whose coefficients are defined by
$\displaystyle{\bar G^1_{11}=H^1_{11},\;G^k_{j1}={g^{ki}\over 2}{\delta g_{ij}
\over\delta t},}$ and
\begin{equation}
\begin{array}{l}\medskip
\displaystyle{L^i_{jk}={g^{im}\over 2}\left({\delta g_{jm}\over\delta x^k}+
{\delta g_{km}\over\delta x^j}-{\delta g_{jk}\over\delta x^m}\right),}\\
\displaystyle{C^{i(1)}_{j(k)}={g^{im}\over 2}\left({\partial g_{jm}\over\partial
y^k}+{\partial g_{km}\over\partial y^j}-{\partial g_{jk}\over\partial y^m}
\right)}.
\end{array}
\end{equation}
By computations, one easily verifies that  $C\Gamma$ satisfies the  conditions
{\it i} and {\it  ii}.

Conversely, let  us consider $\tilde C\Gamma=(\tilde{\bar G}^1_{11},\tilde G^k_{j1},
\tilde L^i_{jk},\tilde C^{i(1)}_{j(k)})$ a h-normal $\Gamma$-linear connection
which satisfies {\it  i} and {\it ii}. It follows directly that
$$
\displaystyle{\tilde{\bar G}^1_{11}=H^1_{11},\;\mbox{and}\;\tilde G^k_{j1}=
{g^{ki}\over 2}{\delta g_{ij}\over\delta t}}.
$$

The condition $g_{ij\vert k}=0$ is equivalent with
$$
{\delta g_{ij}\over\delta x^k}=g_{mj}\tilde L^m_{ik}+g_{im}\tilde L^m_{jk}.
$$
Applying a Christoffel process to the indices $\{i,j,k\}$, we find
$$
\displaystyle{\tilde L^i_{jk}={g^{im}\over 2}\left({\delta g_{jm}\over\delta x^k}+
{\delta g_{km}\over\delta x^j}-{\delta g_{jk}\over\delta x^m}\right)}.
$$

By  analogy, using the relations $C^{i(1)}_{j(k)}=C^{i(1)}_{k(j)}$
and $g_{ij}\vert^{(1)}_{(k)}=0$, following a Christoffel process applied
to the indices $\{i,j,k\}$, we obtain
$$
\displaystyle{\tilde C^{i(1)}_{j(k)}={g^{im}\over 2}\left({\partial g_{jm}\over\partial
y^k}+{\partial g_{km}\over\partial y^j}-{\partial g_{jk}\over\partial y^m}
\right)}.
$$

In conclusion, the uniqueness of the Cartan canonical  connection $C\Gamma$ is
clear. \rule{5pt}{5pt}\medskip\\
\addtocounter{rem}{1}
{\bf Remarks \therem} i) Replacing the canonical nonlinear connection $\Gamma$ by
a general one, the previous theorem holds good.

ii) As a rule, the  Cartan canonical connection of a relativistic rheonomic
Lagrange space $RL^n$ verifies also the properties
\begin{equation}
h_{11/1}=h_{11\vert k}=h_{11}\vert^{(1)}_{(k)}=0\;\mbox{and}\;g_{ij/1}=0.
\end{equation}

iii) Particularly,  the coefficients of the Cartan connection of an autonomous
relativistic rheonomic Lagrange space of electrodynamics (i. e., $g_{ij}(t,x^k,y^k)
=g_{ij}(x^k)$) are the same with those of the Berwald connection, namely,
$C\Gamma=(H^1_{11},0,\gamma^i_{jk},0)$. Note that the Cartan
connection is a $\Gamma$-linear connection, where $\Gamma$ is the canonical
nonlinear connection of the relativistic rheonomic Lagrangian space while the
Berwald connection is a $\Gamma_0$-linear connection, $\Gamma_0$ being the
canonical nonlinear connection  associated to the metric pair $(h_{11},g_{ij})$.
Consequently, the Cartan and Berwald connections are distinct.

iv) The torsion d-tensor {\bf T} of the Cartan canonical connection of a
relativistic  rheonomic Lagrange space is determined by only six local components, because
the properties of the Cartan canonical connection imply the relations
$T^m_{ij}=0$ and $S^{(i)(1)(1)}_{(1)(j)(k)}=0$. At the same time, we point out
that the number of the curvature local d-tensors of the Cartan canonical
connection not reduces. In conclusion, the curvature d-tensor {\bf R} of the
Cartan canonical connection is determined by five effective local d-tensors.
Their expressions was described in Section 1.\medskip\\
\addtocounter{defin}{1}
{\bf Definition \thedefin} The torsion and curvature d-tensors of the Cartan
canonical connection of an $RL^n$ are called the torsion and curvature of
$RL^n$.\medskip

By a direct calculation, we obtain
\begin{th}
i) All torsion d-tensors of an autonomous relativistic rheonomic Lagrange space of
electrodynamics vanish,  except
\begin{equation}
\left\{\begin{array}{l}\medskip
\displaystyle{R^{(m)}_{(1)1j}=-{h_{11}g^{mk}\over 4}\left[H^1_{11}U^{(1)}_{(k)j}+{\partial
U^{(1)}_{(k)j}\over\partial t}\right],}\\
\displaystyle{R^{(m)}_{(1)ij}=r^m_{ijk}y^k+{h_{11}g^{mk}\over 4}\left[U^{(1)}_{(k)i\vert j}+
U^{(1)}_{(k)j\vert i}\right],}
\end{array}\right.
\end{equation}
where $r^m_{ijk}$ are the curvature tensors of the semi-Riemannian metric
$g_{ij}$.\medskip

ii) All curvature d-tensors of an autonomous relativistic rheonomic Lagrange
space of electrodynamics vanish, except $R^l_{ijk}=r^l_{ijk}$.
\end{th}

\section{Relativistic rheonomic Lagrangian electromagnetism}

\subsection{Electromagnetic field}

\setcounter{equation}{0}
\hspace{5mm} Let us consider $RL^n=(J^1(R,M),L)$ a relativistic rheonomic
Lagrange space and $\Gamma=(M^{(i)}_{(1)1},N^{(i)}_{(1)j})$ its canonical
nonlinear connection. At the same time, we denote $C\Gamma=(H^1_{11},G^k_{i1},
L^k_{ij},C^{k(1)}_{i(j)})$ the Cartan canonical connection of $RL^n$.

Using the {\it canonical Liouville d-tensor field} {\bf C}$\displaystyle{=y^i
{\partial\over\partial y^i}}$, we can introduce the {\it deflection
d-tensors}
\begin{equation}
\bar D^{(i)}_{(1)1}=y^i_{/1},\quad
D^{(i)}_{(1)j}=y^i_{\vert j},\quad
d^{(i)(1)}_{(1)(j)}=y^i\vert^{(1)}_{(j)},
\end{equation}
where $"_{/1}"$, $"_{\vert j}"$ and $"\vert^{(1)}_{(j)}"$ are the local
covariant derivatives induced by $C\Gamma$.

By a direct calculation, we find
\begin{prop}
The deflection d-tensors of the rheonomic  Lagrange space $RL^n$ have the
expressions
\begin{equation}\label{dt}
\begin{array}{l}\medskip
\displaystyle{\bar D^{(i)}_{(1)1}={g^{ik}\over 2}{\delta g_{km}\over\delta t}y^m,}
\\\medskip
D^{(i)}_{(1)j}=-N^{(i)}_{(1)j}+L^i_{jm}y^m,
\\
d^{(i)(1)}_{(1)(j)}=\delta^i_j+C^{i(1)}_{m(j)}y^m.
\end{array}
\end{equation}
\end{prop}
\addtocounter{rem}{1}
{\bf Remark \therem} For an autonomous relativistic rheonomic Lagrange space
of electrodynamics
(i. e., $g_{ij}=g_{ij}(x^k)$), the deflection d-tensors reduce to
\begin{equation}
\bar D^{(i)}_{(1)1}=0,\quad
D^{(i)}_{(1)j}=-{1\over 4}g^{ik}h_{11}U^{(1)}_{(k)j},\quad
d^{(i)(1)}_{(1)(j)}=\delta^i_j.
\end{equation}

Using the vertical fundamental metrical d-tensor $G^{(1)(1)}_{(i)(k)}=h^{11}g_{ij}$ of the
relativistic  rheonomic Lagrange space $RL^n$ we construct the {\it metrical
deflection d-tensors},
\begin{equation}
\begin{array}{l}\medskip
\bar D^{(1)}_{(i)1}=G^{(1)(1)}_{(i)(k)}\bar D^{(k)}_{(1)1}=y_{i/1}\\
\medskip
D^{(1)}_{(i)j}=G^{(1)(1)}_{(i)(k)}D^{(k)}_{(1)j}=y_{i\vert j}\\
\medskip
d^{(1)(1)}_{(i)(j)}=G^{(1)(1)}_{(i)(k)}d^{(k)(1)}_{(1)(j)}=y_i\vert^{(1)}_
{(j)},
\end{array}
\end{equation}
where $y_i=G^{(1)(1)}_{(i)(k)}y^k=h^{11}g_{ik}y^k$. Using the
expressions \ref{dt} of the deflection d-tensors, it follows
\begin{prop}
The metrical deflection d-tensors of the relativistic rheonomic Lagrange
space $RL^n$ are given by the formulas
\begin{equation}
\begin{array}{l}\medskip
\displaystyle{\bar D^{(1)}_{(i)1}={h^{11}\over 2}{\delta g_{im}\over\delta t}y^m,}
\\\medskip
D^{(1)}_{(i)j}=h^{11}g_{ik}\left[-N^{(k)}_{(1)j}+L^k_{jm}y^m\right],
\\
d^{(1)(1)}_{(i)(j)}=h^{11}\left[g_{ij}+g_{ik}C^{k(1)}_{m(j)}y^m
\right].
\end{array}
\end{equation}
\end{prop}
\addtocounter{rem}{1}
{\bf Remark \therem} In the particular case of an autonomous relativistic
rheonomic Lagrange space of electrodynamics (i. e., $g_{ij}=g_{ij}(x^k)$),
we have
\begin{equation}
\bar D^{(1)}_{(i)1}=0,\quad
D^{(1)}_{(i)j}=-{1\over 4}U^{(1)}_{(i)j},\quad
d^{(1)(1)}_{(i)(j)}=h^{11}g_{ij}.
\end{equation}

In order to construct the relativistic  rheonomic Lagrangian theory of electromagnetism, we
introduce the following\medskip\\
\addtocounter{defin}{1}
{\bf Definition \thedefin} The distinguished 2-form on $E=J^1(R,M)$
\begin{equation}
F=F^{(1)}_{(i)j}\delta y^i\wedge dx^i+f^{(1)(1)}_{(i)(j)}\delta y^i\wedge
\delta y^j,
\end{equation}
where
\begin{equation}
F^{(1)}_{(i)j}=\displaystyle{{1\over 2}\left[D^{(1)}_{(i)j}-D^{(1)}_{(j)i}
\right]},\quad
f^{(1)(1)}_{(i)(j)}=\displaystyle{{1\over 2}\left[d^{(1)(1)}
_{(i)(j)}-d^{(1)(1)}_{(j)(i)}\right]},
\end{equation}
is called the {\it electromagnetic d-form} of the relativistic rheonomic
Lagrange space $RL^n$.\medskip

Using the above definition, by a direct calculation, we obtain
\begin{prop}
The expressions of the electromagnetic components
\begin{equation}
\left\{\begin{array}{l}\medskip
\displaystyle{F^{(1)}_{(i)j}={h^{11}\over 2}\left[g_{jm}N^{(m)}_{(1)i}-g_{im}N^{(m)}_{(1)j}
+(g_{ik}L^k_{jm}-g_{jk}L^k_{im})y^m\right],}\\
f^{(1)(1)}_{(i)(j)}=0
\end{array}\right.
\end{equation}
hold good.
\end{prop}
\addtocounter{rem}{1}
{\bf Remark \therem} We emphasize that, in the particular case of an autonomous
relativistic rheonomic Lagrange space (i. e. $g_{ij}=g_{ij}(x^k)$), the
electromagnetic local components get the following form
\begin{equation}
\left\{\begin{array}{l}\medskip
\displaystyle{F^{(1)}_{(i)j}={1\over 8}\left[U^{(1)}_{(j)i}-U^{(1)}_{(i)j}\right]}\\
f^{(1)(1)}_{(i)(j)}=0.
\end{array}\right.
\end{equation}

\subsection{Maxwell equations}

\setcounter{equation}{0}
\hspace{5mm} The main result of the electromagnetic relativistic rheonomic Lagrangian
geometry is the following

\begin{th} The electromagnetic local components $F^{(1)}_{(i)j}$ of a
relativistic rheonomic Lagrange space $RL^n=(J^1(R,M),L)$ are governed by
the Maxwell equations
$$
F^{(1)}_{(i)k/1}={1\over 2}{\cal A}_{\{i,k\}}\left\{\bar D^{(1)}_{(i)1\vert k}
+D^{(1)}_{(i)m}T^m_{1k}+d^{(1)(1)}_{(i)(m)}R^{(m)}_{(1)1k}-\left[T^p_{1i\vert k}
+C^{p(1)}_{k(m)}R^{(m)}_{(1)1i}\right] y_p\right\},
$$
$$
\sum_{\{i,j,k\}}F^{(1)}_{(i)j\vert k}=-{1\over 2}\sum_{\{i,j,k\}}C^{(1)(1)(1)}_
{(i)(l)(m)}R^{(m)}_{(1)jk}y^l,\quad\quad\sum_{\{i,j,k\}}F^{(1)}_{(i)j}\vert^{(1)}_
{(k)}=0,
$$
where $y_p=G^{(1)(1)}_{(p)(q)}y^q$ and  $\displaystyle{C^{(1)(1)(1)}_{(i)(l)(m)}
=G^{(1)(1)}_{(l)(q)}C^{q(1)}_{i(m)}={h^{11}\over 2}{\partial^3L\over\partial
y^i\partial y^l\partial y^m}}$.
\end{th}
{\bf Proof.} Firstly, we point out that the Ricci identities \cite{18} applied
to the spatial metrical d-tensor $g_{ij}$ imply that the following curvature
d-tensor identities
$$
R_{mi1k}+R_{im1k}=0,\quad
R_{mijk}+R_{imjk}=0,\quad
P^{\;\;\;\;\;(1)}_{mij(k)}+P^{\;\;\;\;\;(1)}_{imj(k)}=0,
$$
where $R_{mi1k}=g_{ip}R^p_{m1k}$,
$R_{mijk}=g_{ip}R^p_{mjk}$ and $P^{\;\;\;\;\;(1)}_{mij(k)}=g_{ip}
P^{p\;\;\;(1)}_{mj(k)}$, are true.

Now, let us consider the following general deflection d-tensor identities
\cite{18}\medskip

$d_1)\;\;\bar D^{(p)}_{(1)1\vert k}-D^{(p)}_{(1)k/1}=y^mR^p_{m1k}-D^{(p)}_
{(1)m}T^m_{1k}-d^{(p)(1)}_{(1)(m)}R^{(m)}_{(1)1k},$\medskip

$d_2)\;\;D^{(p)}_{(1)j\vert k}-D^{(p)}_{(1)k\vert j}=y^mR^p_{mjk}-d^{(p)(1)}_
{(1)(m)}R^{(m)}_{(1)jk},$\medskip

$d_3)\;\;D^{(p)}_{(1)j}\vert^{(1)}_{(k)}-d^{(p)(1)}_{(1)(k)\vert j}=y^mP^
{p\;\;(1)}_{mj(k)}-D^{(p)}_{(1)m}C^{m(1)}_{j(k)}-d^{(p)(1)}_{(1)(m)}P^{(m)\;
(1)}_{(1)j(k)}$.\medskip\\
Contracting these deflection d-tensor identities by $G^{(1)(1)}_{(i)(p)}$
and using the above curvature d-tensor equalities, we obtain the following
metrical deflection d-tensors identities:\medskip

$d^\prime_1)\;\;\bar D^{(1)}_{(i)1\vert k}-D^{(1)}_{(i)k/1}=-y_mR^m_{i1k}-
D^{(1)}_{(i)m}T^m_{1k}-d^{(1)(\mu)}_{(i)(m)}R^{(m)}_{(1)1k},$\medskip

$d^\prime_2)\;\;D^{(1)}_{(i)j\vert k}-D^{(1)}_{(i)k\vert j}=-y_mR^m_{ijk}-
d^{(1)(\mu)}_{(i)(m)}R^{(m)}_{(1)jk},
$\medskip

$d^\prime_3)\;\;D^{(1)}_{(i)j}\vert^{(1)}_{(k)}-d^{(1)(1)}_{(i)(k)\vert j}=
-y_mP^{m\;\;(1)}_{ij(k)}-D^{(1)}_{(i)m}C^{m(1)}_{j(k)}-d^{(1)(1)}_{(i)(m)}
P^{(m)\;(1)}_{(1)j(k)}$.\medskip

At the same time, we recall that the following Bianchi identities \cite{15}
\medskip

$b_1)\;\;{\cal A}_{\{j,k\}}\left\{R^l_{j1k}+T^l_{1j\vert k}+C^{l(1)}_{k(m)}
R^{(m)}_{(1)1j}\right\}=0$,\medskip

$b_2)\;\;\sum_{\{i,j,k\}}\left\{R^l_{ijk}-C^{l(1)}_{k(m)}R^{(m)}_{(1)ij}
\right\}=0$,\medskip

$b_3)\;\;{\cal A}_{\{j,k\}}\left\{P^{l\;\;(1)}_{jk(p)}+C^{l(1)}_{j(p)\vert k}
+C^{l(1)}_{k(m)}P^{(m)\;\;(1)}_{(1)j(p)}\right\}=0$,
\medskip\\
where ${\cal A}_{\{j,k\}}$ means alternate sum and $\sum_{\{i,j,k\}}$ means
cyclic sum, hold good.

In order to obtain the first Maxwell identity, we permute $i$ and $k$ in
$d^\prime_1$ and we subtract the new identity from the initial one. Finally,
using the Bianchi identity $b_1$, we obtain what we were looking for.

Doing a cyclic sum by the indices $\{i,j,k\}$ in $d^\prime_2$ and using the
Bianchi identity $b_2$, it follows the second Maxwell equation.

Applying a Christoffel process to the indices $\{i,j,k\}$ id $d^\prime_3$
and combining with the Bianchi identity $b_3$ and the relation
$P^{(m)\;\;(1)}_{(1)j(p)}=P^{(m)\;\;(1)}_{(1)p(j)}$,
we get a new identity. The cyclic sum by the indices $\{i,j,k\}$ applied to
this last identity implies the third Maxwell equation.
\rule{5pt}{5pt}\medskip\\
\addtocounter{rem}{1}
{\bf Remark \therem} In the case of an autonomous relativistic rheonomic Lagrange space of
electrodynamics (i. e., $g_{ij}=g_{ij}(x^k)$), the Maxwell equations take the
simple form
\begin{equation}
\hspace*{5mm}
F^{(1)}_{(i)k/1}={1\over 2}{\cal A}_{\{i,k\}}h^{11}g_{im}R^{(m)}_{(1)1k},
\quad\sum_{\{i,j,k\}}F^{(1)}_{(i)j\vert k}=0,\quad
\sum_{\{i,j,k\}}F^{(1)}_{(i)j}\vert^{(1)}_{(k)}=0.
\end{equation}

\section{Relativistic rheonomic Lagrangian gravitational theory}

\subsection{Gravitational field}
\setcounter{equation}{0}
\hspace{5mm} Let $h=(h_{11})$ be a fixed semi-Riemannian metric on
the temporal manifold $R$ and $\Gamma=(M^{(i)}_{(1)1},N^{(i)}_{(1)j})$ a
fixed nonlinear connection on the 1-jet space $J^1(R,M)$. In order to
develope a relativistic rheonomic Lagrange theory of gravitational field
on $J^1(R,M)$, we introduce the following\medskip\\
\addtocounter{defin}{1}
{\bf Definition \thedefin} From physical point of view, an adapted metrical d-tensor
$G$ on $J^1(R,M)$, expressed locally by
$$
G=h_{11}dt\otimes dt+g_{ij}dx^i\otimes dx^j+h^{11}g_{ij}\delta y^i\otimes
\delta y^j,
$$
where $g_{ij}=g_{ij}(t,x^k,y^k)$ is a d-tensor on $J^1(R,M)$, symmetric,
of rank $n=\dim M$ and having a constant signature on $E$, is called a
{\it gravitational $h$-potential} on $E$.

Now, taking $RL^n=(J^1(R,M),L)$ a relativistic rheonomic Lagrange space, via
its vertical fundamental metrical d-tensor
$$
G^{(1)(1)}_{(i)(j)}={1\over 2}{\partial^2L\over\partial y^i\partial y^j}=
h^{11}(t)g_{ij}(t,x^k,y^k),
$$
and its canonical nonlinear connection $\Gamma=(M^{(i)}_{(1)1},N^{(i)}_{(1)j})$,
one induces a natural gravitational $h$-potential on $J^1(R,M)$, setting
\begin{equation}
G=h_{11}dt\otimes dt+g_{ij}dx^i\otimes dx^j+h^{11}g_{ij}\delta y^i\otimes
\delta y^j.
\end{equation}

\subsection{Einstein equations and conservation  laws}

\setcounter{equation}{0}
\hspace{5mm} Let us consider $C\Gamma=(H^1_{11},G^k_{j1},L^i_{jk},C^{i(1)}_{j(k)})$ the
Cartan canonical connection of the relativistic rheonomic Lagrange space
$RL^n$.

We postulate that the Einstein equations which govern the gravitational $h$-potential
$G$ of the relativistic rheonomic Lagrange space $RL^n$ are the Einstein equations attached
to the Cartan canonical connection of $RL^n$ and the adapted metric $G$ on
$J^1(R,M)$, that is,
\begin{equation}
Ric(C\Gamma)-{Sc(C\Gamma)\over 2}G={\cal K}{\cal T},
\end{equation}
where $Ric(C\Gamma)$ represents the Ricci d-tensor of the Cartan connection,
$Sc(C\Gamma)$ is its scalar curvature, ${\cal K}$ is the Einstein constant
and ${\cal T}$ is an intrinsec tensor of matter which is called  the
{\it stress-energy} d-tensor.

In the adapted basis $(X_A)=\displaystyle{\left({\delta\over\delta t},
{\delta\over\delta x^i},{\partial\over\partial y^i}\right)}$, the curvature
d-tensor {\bf R} of the Cartan connection is expressed locally by {\bf R}$(X_
C,X_B)X_A=R^D_{ABC}X_D$. Hence, it follows that we have $R_{AB}=Ric(C\Gamma)
(X_A,X_B)=R^D_{ABD}$ and $Sc(C\Gamma)=G^{AB}R_{AB}$, where
\begin{equation}
G^{AB}=\left\{\begin{array}{ll}\medskip
h_{11},&\mbox{for}\;\;A=1,\;B=1\\\medskip
g^{ij},&\mbox{for}\;\;A=i,\;B=j\\\medskip
h_{11}g^{ij},&\mbox{for}\;\;A={(i)\atop(1)},\;B={(j)\atop(1)}\\
0,&\mbox{otherwise}.
\end{array}\right.
\end{equation}

Taking into account the expressions of the local curvature d-tensors of the
Cartan connection and the form of the vertical fundamental metrical d-tensor
$G^{AB}$, we deduce
\begin{prop}
The Ricci d-tensor of the Cartan canonical  connection of $RL^n$ is determined
by the following six  effective local Ricci d-tensors,
$$
R_{11}\stackrel{\mbox{not}}{=}H_{11}=H^1_{111}=0,\quad
R_{i1}=R^m_{i1m},\quad R_{ij}=R^m_{ijm},\quad
R^{\;(1)}_{i(j)}\stackrel{\mbox{not}}{=}P^{\;(1)}_{i(j)}=-P^{m\;\;(1)}_{im(j)},
$$
$$
R^{(1)}_{(i)1}\stackrel{\mbox{not}}{=}P^{(1)}_{(i)1}=P^{m\;\;(1)}_{i1(m)},\quad
R^{(1)}_{(i)j}\stackrel{\mbox{not}}{=}P^{(1)}_{(i)j}=P^{m\;\;(1)}_{ij(m)},\quad
R^{(1)(1)}_{(i)(j)}\stackrel{\mbox{not}}{=}S^{(1)(1)}_{(i)(j)}=S^{m(1)(1)}_{i(j)(m)}.
$$
\end{prop}

Consequently, denoting $H=h^{11}H_{11},\;\;R=g^{ij}R_{ij}$ and $S=h_{11}g^{ij}
S^{(1)(1)}_{(i)(j)}$, we obtain
\begin{prop}
The scalar curvature of the Cartan canonical connection of $RL^n$ has the
expression
\begin{equation}
Sc(C)=H+R+S=R+S.
\end{equation}
\end{prop}
\addtocounter{rem}{1}
{\bf Remark \therem} In the particular case of an autonomous relativistic
rheonomic Lagrange space of electrodynamics (i. e., $g_{ij}=g_{ij}(x^k)$),
all Ricci d-tensors vanish, except $R_{ij}=r_{ij}$, where $r_{ij}$ are the
Ricci tensors associated to the semi-Riemannian metric $g_{ij}$. It follows
that the scalar curvatures of a such space are $H=0,\;R=r,\;S=0$, where $H$
and $r$ are the scalar curvatures of the semi-Riemannian metrics $h_{11}$
and $g_{ij}$. \medskip

Using the above results, we can establish the following

\begin{th}
The Einstein equations which govern the gravitational $h$-potential
$G$ induced by the Kronecker $h$-regular Lagrangian function of a relativistic
rheonomic Lagrange space $RL^n$, have the form
$$
\left\{\begin{array}{l}\medskip
\displaystyle{-{R+S\over 2}h_{11}={\cal K}{\cal T}_{11}}\\\medskip
\displaystyle{R_{ij}-{R+S\over 2}g_{ij}={\cal K}{\cal T}_{ij}}\\\medskip
\displaystyle{S^{(1)(1)}_{(i)(j)}-{R+S\over 2}h^{11}g_{ij}={\cal K}{\cal T}
^{(1)(1)}_{(i)(j)}},
\end{array}\right.\leqno{(E_1)}
$$
$$
\left\{\begin{array}{lll}\medskip
0={\cal T}_{1i},&R_{i1}={\cal K}{\cal T}_{i1},&
P^{(1)}_{(i)1}={\cal K}{\cal T}^{(1)}_{(i)1}\\
0={\cal T}^{\;(1)}_{1(i)},&P^{\;(1)}_{i(j)}={\cal K}{\cal T}^{\;(1)}_{i(j)},&
P^{(1)}_{(i)j}={\cal K}{\cal T}^{(1)}_{(i)j},
\end{array}\right.\leqno{(E_2)}
$$
where ${\cal T}_{AB},\;A,B\in\{1,i,{(1)\atop(i)}\}$ are the adapted
local components of the stress-energy d-tensor ${\cal T}$.
\end{th}
\addtocounter{rem}{1}
{\bf Remark \therem} i) Note that, in order to have the compatibility of the Einstein equations, it
is necessary that the certain adapted local components of the stress-energy
d-tensor vanish {\it "a priori"}.

ii) In the particular case of an autonomous relativistic
rheonomic Lagrange space of electrodynamics (i. e., $g_{ij}=g_{ij}(x^k)$), using
preceding notations, the following Einstein equations of gravitational field,
$$
\left\{\begin{array}{l}\medskip
\displaystyle{r_{ij}-{r\over 2}g_{ij}={\cal K}{\cal T}_{ij}}\\
\displaystyle{-{r\over 2}h^{11}g_{ij}={\cal K}{\cal T}^{(1)(1)}_{(i)(j)}},
\end{array}\right.\leqno{(E_1)}
$$
$$
\left\{\begin{array}{lll}\medskip
0={\cal T}_{1i},&0={\cal T}_{i1},&
0={\cal T}^{(1)}_{(i)1}\\
0={\cal T}^{\;(1)}_{1(i)},&0={\cal T}^{\;(1)}_{i(j)},&
0={\cal T}^{(1)}_{(i)j},
\end{array}\right.\leqno{(E_2)}
$$
hold good.\medskip

It is well known that, from physical point of view, the
stress-energy d-tensor ${\cal T}$ must verifies the local {\it conservation
laws} ${\cal T}^B_{A\vert B}=0,\;\;\forall\;A\in\{1,i,{(1)\atop(i)}\}$, where
${\cal T}^B_A=G^{BD}{\cal T}_{DA}$.

\begin{th}
In the relativistic rheonomic Lagrangian geometry, the {\it conservation laws} of
the Einstein equations are
\begin{equation}
\left\{\begin{array}{l}\medskip
\displaystyle{\left[{R+S\over 2}\right]_{/1}=R^m_{1\vert m}-P^{(m)}_{(1)1}
\vert^{(1)}_{(m)}}\\\medskip
\displaystyle{\left[R^m_j-{R+S\over 2}\delta^m_j\right]_{\vert m}=
-P^{(m)}_{(1)j}\vert^{(1)}_{(m)}}\\\medskip
\displaystyle{\left[S^{(m)(1)}_{(1)(j)}-{R+S\over 2}\delta^m_j\right]\left\vert
^{(1)}_{(m)}\right.=-P^{m(1)}_{\;\;\;(j)\vert m}},
\end{array}\right.
\end{equation}
where\medskip
$
\;R^i_1=g^{im}R_{m1},\;\; P^{(i)}_{(1)1}=h_{11}g^{im}P^{(1)}_{(m)1},\;\;
R^i_j=g^{im}R_{mj},\;\;P^{(i)}_{(1)j}=h_{11}g^{im}P^{(1)}_{(m)j},
$
$
P^{i(1)}_{\;\;(j)}=g^{im}P^{\;\;(1)}_{m(j)}\;\mbox{and}\;
S^{(i)(1)}_{(1)(j)}=h_{11}g^{im}S^{(1)(1)}_{(m)(j)}.
$
\end{th}
{\bf Acknowledgements.} The author thank to the reviewers of the Duke Mathematical
Journal for their valuable comments upon a previous version of this paper.

\begin{center}
University POLITEHNICA of Bucharest\\
Department of Mathematics I\\
Splaiul Independentei 313\\
77206 Bucharest, Romania\\
e-mail: mircea@mathem.pub.ro\\
\end{center}
\end{document}